\documentclass{amsart}

\usepackage{amsmath,amsfonts,amssymb,amsthm,graphicx,cite,mathrsfs}

\newtheorem{theorem}{Theorem}[section]
\newtheorem{lemma}{Lemma}[section]
\newtheorem{proposition}{Proposition}[section]

\theoremstyle{remark}
\newtheorem{remark}{Remark}[section]

\begin{document}

\title[Multivariable Bessel polynomials]{An orthogonality relation for multivariable Bessel polynomials}

\author{Martin Halln\"as} \address{SISSA, Via Beirut 2-4, 34014 Trieste TS, Italy}
\email{hallnas@sissa.it}
\thanks{The author was supported by the European Union through the FP6 Marie Curie RTN ENIGMA (Contract number MRTN-CT-200405652).}

\date{\today}
\subjclass[2000]{Primary 33C50}
\keywords{Bessel polynomials, multivariable orthogonal polynomials}

\begin{abstract}
In a recent paper we introduced a multivariable generalisation of the Bessel polynomials, depending on one extra parameter, and related to the so-called hyperbolic Sutherland model with external Morse potential. In this paper we obtain a corresponding multivariable generalisation of a well-known orthogonality relation for the (one-variable) Bessel polynomials due to Krall and Frink.
\end{abstract}

\maketitle

\section{Introduction}
Up to normalisation, the Bessel polynomials $y_m(x;a)$, which depend on a (complex) parameter $a$ and a non-negative integer $m$, can be defined as the polynomials of degree $m$ such that they are eigenfunctions of the ordinary differential operator
\begin{equation*}
	d^B = x^2\frac{\partial^2}{\partial x^2} + (ax + 2)\frac{\partial}{\partial x}.
\end{equation*}
As such, they first appeared in a paper by Bochner \cite{Boc29}, although the name Bessel polynomials was introduced later by Krall and Frink \cite{KF49}. We recall that a second parameter $b$ can be introduced by a change of coordinates to $y = bx$. As first shown by Hahn \cite{Hah35}, there exist no interval on the real line and corresponding non-negative weight function with respect to which the Bessel polynomials are orthogonal. However, Krall \cite{Kra41} demonstrated that the Bessel polynomials are orthogonal in a more general sense. In fact, they satisfy the orthogonality relation
\begin{equation}\label{oneVarOrthogonalityRel}
	\int_{|x|=1}y_n(x)y_m(x)\rho(x)dx = 0,\quad n\neq m,
\end{equation}
where
\begin{equation*}
	\rho(x;a) = \frac{1}{2\pi i}\sum_{n=0}^\infty \frac{\Gamma(a)}{\Gamma(a+n-1)}\left(-\frac{2}{x}\right)^n.
\end{equation*}
The result was first stated in this form by Krall and Frink \cite{KF49}. For a general account of the theory of Bessel polynomials, and an excellent guide to the extensive literature on the topic (up to 1978), see the book by Grosswald \cite{Gro78}.

In a recent paper we introduced and studied a multivariable generalisation of the Bessel polynomials as eigenfunctions of the partial differential operator
\begin{equation*}
	D^B = \sum_{i=1}^n\frac{\partial^2}{\partial x_i^2} + \sum_{i=1}^n(ax_i + 2)\frac{\partial}{\partial x_i} + 2\kappa\sum_{i\neq j}\frac{x_i^2}{x_i-x_j}\frac{\partial}{\partial x_i},
\end{equation*}
where $\kappa$ is a (complex) parameter; see Halln\"as \cite{Hal08}. In this latter paper we emphasised the relation with a quantum many-body system of so-called Calogero-Moser-Sutherland type. This lead to certain orthogonality results with respect to the inner product in a particular Hilbert space of square-integrable functions, to be discussed further in Section \ref{proofThm21}. However, these orthogonality results are valid only for a non-generic set of parameter values, and then involve only a finite set of multivariable Bessel polynomials. The main purpose of this paper is to deduce a multivariable generalisation of the orthogonality relation \eqref{oneVarOrthogonalityRel}, valid for generic parameter values, and involving all multivariable Bessel polynomials.

We conclude this introduction by mentioning that Duran \cite{Dur93}, Evans et al.~\cite{EEKL93}, as well as Kwon et al.~\cite{KKH92}, have obtained real, but necessarily not non-negative, weight functions with respect to which the (one-variable) Bessel polynomials are orthogonal on $\lbrack 0,\infty)$. The problem of obtaining corresponding multivariable weight functions remains open.

\section{Orthogonality relation and normalisation factors}\label{statementsSection}
In this section we state our main result: the orthogonality relation for our multivariable Bessel polynomials. In order to do so we first briefly review some basic results on symmetric functions and Jack polynomials, recall the definition of our multivariable Bessel polynomials, and state a result on the nature of their dependence on the two parameters $a$ and $\kappa$.

\subsection{Symmetric functions and Jack polynomials}\label{symFuncsJackPolsSec}
A partition $\lambda = (\lambda_1,\ldots,\lambda_n)$ is a sequence of non-negative integers such that $\lambda_1\geq\cdots\geq\lambda_n$, and its weight is the sum $|\lambda| = \lambda_1+\cdots+\lambda_n$ of its parts. On the set of partitions of a given fixed weight the so-called dominance order is defined by
\begin{equation*}
	\mu\leq\lambda\Leftrightarrow \mu_1+\cdots+\mu_i\leq\lambda_1+\cdots+\lambda_i,\quad \forall i=1,\ldots,n.
\end{equation*}
The monomial symmetric polynomials $m_\lambda$ in the variables $x = (x_1,\ldots,x_n)$ are given by
\begin{equation*}
	m_\lambda(x) = \sum_\alpha x_1^{\alpha_1}\cdots x_n^{\alpha_n},
\end{equation*}
where the sum is over all distinct permutations $\alpha$ of $\lambda$. These monomials form a linear basis for the algebra of symmetric polynomials in $n$ variables with complex coefficients, from hereon denoted $\Lambda_n$. As long as the parameter $\kappa$ is not a negative rational number, there exists for each partition $\lambda = (\lambda_1,\ldots,\lambda_n)$ a unique symmetric polynomial $P_\lambda$ in the variables $x$, referred to as a Jack polynomial, of the form
\begin{equation*}
	P_\lambda = m_\lambda + \sum_{\mu<\lambda}u_{\lambda\mu}m_\mu
\end{equation*}
for some coefficients $u_{\lambda\mu}$ such that it is an eigenfunction of the partial differential operator
\begin{equation*}
	D = \sum_{i=1}^n x_i^2\frac{\partial^2}{\partial x_i^2} + 2\kappa\sum_{i\neq j}\frac{x_i^2}{x_i-x_j}\frac{\partial}{\partial x_i};
\end{equation*}
see e.g.~Section 3 in Stanley \cite{Sta89}. We remark that the parameter $\kappa$ is the inverse of the parameter $\alpha$ used by Stanley. We let $T^n = \lbrace x\in\mathbb{C}^n: |x_i| = 1, \forall i =1,\ldots,n\rbrace$ be the n-dimensional torus. In addition, we will use the notation $x^{-1} = (1/x_1,\ldots,1/x_n)$, $dx = dx_1\cdots dx_n$ and $dx/x = dx_1/x_1\cdots dx_n/x_n$. If $\Re(\kappa)\geq 0$ then the Jack polynomials are orthogonal polynomials with respect to the bilinear form on $\Lambda_n$ defined by
\begin{equation*}
	\langle f,g\rangle^\prime_n = \frac{1}{(2\pi i)^nn!}\int_{T^n}f(x)g(x^{-1})\prod_{i\neq j}\left(1 - \frac{x_i}{x_j}\right)^\kappa \frac{dx}{x}
\end{equation*}
for any $f,g\in\Lambda_n$, and the corresponding normalisation factors are given by
\begin{equation*}
	\langle P_\lambda,P_\lambda\rangle^\prime_n = \prod_{1\leq i<j\leq n}\frac{\Gamma(\kappa(j-i+1)+\lambda_i-\lambda_j)\Gamma(\kappa(j-i-1)+\lambda_i-\lambda_j+1)}{\Gamma(\kappa(j-i)+\lambda_i-\lambda_j)\Gamma(\kappa(j-i)+\lambda_i-\lambda_j+1)};
\end{equation*}
see e.g.~Section VI.10 in Macdonald \cite{Mac95}. Here, as well as in the discussions below, we use the convention that empty products are equal to one. For a detailed exposition of the theory of symmetric functions and Jack polynomials we refer to the book by Macdonald \cite{Mac95}.

\subsection{The multivariable Bessel polynomials}
In analogy with the Jack polynomials, the multivariable Bessel polynomials can now be defined by the following fact (see Section 3 in Halln\"as \cite{Hal08}): given a partition $\lambda = (\lambda_1,\ldots,\lambda_n)$ there exists for generic values of the parameters $a$ and $\kappa$ a unique symmetric polynomial $Y_\lambda$ of the form
\begin{equation*}
	Y_\lambda = P_\lambda + \sum_{\mu\subset\lambda}u_{\lambda\mu}P_\mu
\end{equation*}
for some coefficients $u_{\lambda\mu}$ such that it is an eigenfunction of $D^B$. We remark that it is rather easy to establish the existence of such a polynomial for generic parameter values, but to obtain a condition on the parameters which is both necessary and sufficient is much harder. In Halln\"as (see Proposition 3.1) we obtained a sufficient condition which, however, was rather far from being necessary. To state our main result we will need such a condition which is much closer to being necessary. We will resolve this problem by proving the following:

\begin{proposition}\label{rationalityProp}
The multivariable Bessel polynomials $Y_\lambda(a,\kappa)$ are rational functions of $a$ and $\kappa$ with poles only along hyperplanes given by $\kappa $ a negative rational number or $a -1 + \kappa i$ a non-positive integer for some $i = 0,1,\ldots,2(n-1)$.
\end{proposition}

\subsection{Statement of main result}
We will show that the multivariable Bessel polynomials are orthogonal polynomials on the torus $T^n$ with respect to the weight function
\begin{equation*}
	W(x;a,\kappa) = \mathscr{W}(x;a,\kappa)\prod_{i\neq j}\left(1 - \frac{x_i}{x_j}\right)^\kappa,
\end{equation*}
where
\begin{multline}\label{WDef}
	\mathscr{W}(x;a,\kappa) = \frac{1}{(2\pi i)^n n!}\sum_\lambda \Bigg(\prod_{i<j}\frac{\Gamma(\kappa(j-i)+\lambda_i-\lambda_j+1)}{\Gamma(\kappa(j-i-1)+\lambda_i-\lambda_j+1)}\\ \times\prod_{i=1}^n \frac{\Gamma(\kappa i)\Gamma(a+\kappa(2n-i-1))}{\Gamma(\kappa)\Gamma(a-1+\kappa(2n-i-1)+\lambda_i)}\Bigg) P_\lambda(-2x^{-1};\kappa).
\end{multline}
A priori, it is not clear for which $x\in\mathbb{C}^n$ the right-hand side of \eqref{WDef} converges, or indeed if such $x$ exist. This problem is settled by the following:

\begin{lemma}\label{convergenceLemma}
Assume that $\kappa$ is not a negative rational number and that $a-1+\kappa(2n-i-1)$ is not a negative integer for any $i=1,\ldots,n$. Let $N$ be a neighbourhood of $0\in\mathbb{C}^n$. Then the right-hand side of \eqref{WDef} converges uniformly and absolutely on $\mathbb{C}^n\setminus N$.
\end{lemma}

We let $\lbrack\alpha\rbrack_m$ denote the Pochhammer symbol, defined by $\lbrack\alpha\rbrack_0 = 1$ and $\lbrack\alpha\rbrack_m = \alpha(\alpha+1)\cdots(\alpha+m-1)$ for all positive integers $m$. In addition, we will use the notation $\lbrack\alpha_1,\ldots,\alpha_p\rbrack_m = \lbrack\alpha_1\rbrack_m\cdots\lbrack\alpha_p\rbrack_m$. We are now ready to state our main result.

\begin{theorem}\label{mainThm}
Assume that $\Re(\kappa)\geq 0$ and that $a-1+\kappa i\neq 0,-1,\ldots$ for all $i = 0,\ldots,2(n-1)$. Let $\lambda = (\lambda_1,\ldots,\lambda_n)$ and $\mu = (\mu_1,\ldots,\mu_n)$ be two partitions such that $\lambda\neq \mu$. Then
\begin{equation}\label{orthogonalityRel}
	\int_{T^n}Y_\lambda(x)Y_\mu(x)W(x)dx = 0.
\end{equation}
Moreover,
\begin{multline}\label{normFactorFormula}
	\int_{T^n}Y_\lambda^2(x)W(x)dx\\ = (-1)^{|\lambda|+n}2^{2|\lambda|+n} \prod_{1\leq i<j\leq n}\Bigg(\frac{\lbrack\kappa(j-i+1), \kappa(j-i-1)+1\rbrack_{\lambda_i-\lambda_j}}{\lbrack\kappa(j-i), \kappa(j-i)+1\rbrack_{\lambda_i-\lambda_j}}\\ \shoveright{\times\frac{\lbrack a+\kappa(2n-i-j+1)-1, a+\kappa(2n-i-j-1)\rbrack_{\lambda_i+\lambda_j}}{\lbrack a+\kappa(2n-i-j)-1, a+\kappa(2n-i-j)\rbrack_{\lambda_i+\lambda_j}}\Bigg)}\\ \times\prod_{i=1}^n \frac{\lbrack a+\kappa(n-i)-1, \kappa(n-i)+1\rbrack_{\lambda_i}}{\lbrack a+2\kappa(n-i)-1, a+2\kappa(n-i)\rbrack_{2\lambda_i}}.
\end{multline}
\end{theorem}

\begin{remark}
Suppose that $n = 1$. Then, as is readily inferred from their definitions, $P_{(m)}(x_1) = x_1^m$ and $Y_{(m)}(x_1) = C_my_m(x_1)$ for some constant $C_m$, depending on the particular normalisation chosen for the $y_m$. We recall that Krall and Frink used the normalisation in which the constant term of  each $y_m$ is equal to one. In this case, $C_m = 2^m \lbrack a-1\rbrack_m/\lbrack a-1\rbrack_{2m}$; see Section 6 in Halln\"as \cite{Hal08} for the precise relation between this and the monic normalisation in the multivariable case. It is thus clear that $W(x_1) = \rho(x_1)$, and that we recover Krall and Frink's orthogonality relation \eqref{oneVarOrthogonalityRel} for the Bessel polynomials $y_m$ by setting $n = 1$ in Theorem \ref{mainThm}.
\end{remark}

\section{Proofs}
In this section we prove the results stated in Section \ref{statementsSection}.

\subsection{Proof of Proposition \ref{rationalityProp}}
In a first step we will deduce a closed expression for the coefficients in the expansion of the multivariable Bessel polynomials in Jack polynomials. Our starting point will be such an expression for the Jacobi polynomials associated with the root system $BC_n$, obtained by Okounkov and Olshanski \cite{OO06} (see Proposition 2.3). The reason is that the multivariable Bessel polynomials can be obtained from these $BC_n$ Jacobi polynomials via a limit transition; see Proposition 4.1 in Halln\"as. By applying this limit transition to Okounkov and Olshanski's result we will obtain the desired expression for the multivariable Bessel polynomials. It will then be a straightforward exercise to infer Proposition \ref{rationalityProp} from this expression.

We let $t(z) = (4\sinh^2 z_1/2,\ldots,4\sinh^2 z_n/2)$. For each partition $\lambda = (\lambda_1,\ldots,\lambda_n)$ and generic (complex) values of the parameters $(k_1,k_2,k_3)$ a $BC_n$ Jacobi polynomial $P^{BC}_\lambda$ can be defined as the unique eigenfunction of the partial differential operator
\begin{multline*}
	D^{BC} = \sum_{i=1}^n\frac{\partial^2}{\partial z_i^2} + \sum_{i=1}^n\left(k_1\coth\frac{1}{2}z_i + 2k_2\coth z_i\right)\frac{\partial}{\partial z_i}\\ + k_3\sum_{i<j}\left(\coth\frac{1}{2}(z_i - z_j)\left(\frac{\partial}{\partial z_i} - \frac{\partial}{\partial z_j}\right) + \coth\frac{1}{2}(z_i + z_j)\left(\frac{\partial}{\partial z_i} + \frac{\partial}{\partial z_j}\right)\right)
\end{multline*}
such that
\begin{equation*}
	P^{BC}_\lambda(z) = P_\lambda(t(z)) + \sum_{\mu\subset\lambda}u_{\lambda\mu}P_\mu(t(z))
\end{equation*}
for some coefficients $u_{\lambda\mu}$; see e.g.~Beerends and Opdam \cite{BO93} and references therein. An important ingredient in Okounkov and Olshanski's result is the so-called $BC_n$ interpolation polynomials. We proceed to recall a 'combinatorial' formula for these polynomials which we shall make use of at a later stage in the proof. To this end we let $\mu = (\mu_1,\ldots,\mu_n)$ be a partition, and recall that a reverse column-strict tableux $T$ of shape $\mu$ can be defined as a function that assigns to each box $(i,j)\in\mu$ an integer $T(i,j)\in\lbrace 1,\ldots,n\rbrace$ such that $T(i,j)$ decreases strictly down each column and weakly from left to right along each row. Suppose that $\kappa$ is not a negative rational number. The $BC_n$ interpolation polynomial $I_\mu$ can then be defined by a formula of the form
\begin{multline}\label{interpolationPols}
	I_\mu(x;\kappa,h) = \sum_T \psi_T(\kappa)\prod_{(i,j)\in\mu}\big((x_{T(i,j)}+h-\kappa T(i,j))^2\\ -(j-1-\kappa(i-1)+h-\kappa T(i,j))^2\big),
\end{multline}
where the sum is over all reverse column-strict tableux $T$ of shape $\mu$; see Proposition 2.2 in Okounkov and Olshanski \cite{OO06}. The weights $\psi_T(\kappa)$ are rational functions in $\kappa$ with poles located only at negative rational numbers. An explicit formula for these weights can be found in Section VI.10 of Macdonald \cite{Mac95} (see Equations (10.10)--(10.12)). Since we will not make use of this formula we refrain from stating it here. The results of Proposition 2.3 in Okounkov and Olshanski \cite{OO06} can now be formulated as follows:

\begin{theorem}[Okounkov and Olshanski]
Let $\sigma = (k_1 + 2k_2)/2$. Then
\begin{equation}\label{JacobiExp}
	\frac{P^{BC}_\lambda(z;k_1,k_2,k_3)}{P^{BC}_\lambda(0^n;k_1,k_2,k_3)}\\ = \sum_{\mu\subseteq\lambda}\frac{I_\mu(\lambda;k_3,\sigma+k_3 n)P_\mu(t(z);k_3)}{I_\mu(\mu;k_3,\sigma+k_3 n)P^{BC}_\mu(0^n;k_1,k_2,k_3)}.
\end{equation}
\end{theorem}

\begin{remark}
Okounkov and Olshanski referred to this formula for the $BC_n$ Jacobi polynomials as a binomial formula.
\end{remark}

At this point we set $(k_1,k_2,k_3) = ((a-1+2e^\epsilon)/2,(a-1-2e^\epsilon)/4,\kappa)$ with $\epsilon>0$. In particular, we then have $\sigma = (a-1)/2$. We have, furthermore, the limit transitions
\begin{equation*}
	\lim_{\epsilon\rightarrow\infty} e^{-\epsilon|\lambda|}P^{BC}_\lambda(z_1+\epsilon,\ldots,z_n+\epsilon;k_1,k_2,k_3) = Y_\lambda(e^{z_1},\ldots,e^{z_n};a,
\kappa)
\end{equation*}
and
\begin{equation*}
	\lim_{\epsilon\rightarrow\infty} e^{-\epsilon|\lambda|}P^{BC}_\lambda(0^n;k_1,k_2,k_3) = Y_\lambda(0^n;a,\kappa);
\end{equation*}
see Proposition 4.1 and the proof of Proposition 6.1 in Halln\"as \cite{Hal08}, respectively. We observe that $\lim_{\epsilon\rightarrow\infty} 4e^{-\epsilon}\sinh^2(z+\epsilon)/2 = e^z$. By taking the limit $\epsilon\rightarrow\infty$ in \eqref{JacobiExp} we thus obtain the following:

\begin{proposition}
The multivariable Bessel polynomials have the series expansion
\begin{equation*}
	\frac{Y_\lambda(x;a,\kappa)}{Y_\lambda(0^n;a,\kappa)} = \sum_{\mu\subseteq\lambda}\frac{I_\mu(\lambda;\kappa,(a-1)/2+\kappa n)P_\mu(x;\kappa)}{I_\mu(\mu;\kappa,(a-1)/2+\kappa n)Y_\mu(0^n;a,\kappa)}.
\end{equation*}
\end{proposition}

We recall that the Jack polynomials $P_\lambda(\kappa)$ are rational functions of $\kappa$ with poles only at negative rational numbers; see e.g.~Equation (10.13) and the preceding discussion in Chapter VI of Macdonald \cite{Mac95}. We observe that \eqref{interpolationPols} implies that this fact holds true also for $I_\mu(\lambda;\kappa,(a-1)/2+\kappa n)$. As discussed by Okounkov and Olshanski \cite{OO06}, the $BC_n$ interpolation polynomials $I_\mu$ are normalised such that
\begin{multline*}
	I_\mu(\mu;\kappa,(a-1)/2+\kappa n)\\ = \prod_{(i,j)\in\mu}\big(1+\mu_i-j+\kappa(\mu^\prime_j-i)\big)\big(a-2+\mu_i+j+\kappa(2n-\mu^\prime_j-i)\big).
\end{multline*}
It is clear that a factor $1+\mu_i-j+\kappa(\mu^\prime_j-i) = 0$ only if $\kappa$ is a negative rational number. Since $\mu_i+j-1> 0$ and $2\leq\mu^\prime_j+i\leq 2n$, it is also clear that a factor $a-2+\mu_i+j+\kappa(2n-\mu^\prime_j-i) = 0$ only if $a-1+\kappa i$ is a negative integer for some $i = 0,\ldots,2(n-1)$. Furthermore, the specialisation of the multivariable Bessel polynomials at $x = (0^n)$ is given by
\begin{equation*}
\begin{split}
	Y_\lambda(0^n;a,\kappa) &= 2^{|\lambda|}\prod_{i<j}\frac{\lbrack \kappa(j-i+1)\rbrack_{\lambda_i-\lambda_j}\lbrack \kappa(2n-i-j+1)\rbrack_{\lambda_i+\lambda_j}}{\lbrack \kappa(j-i)\rbrack_{\lambda_i-\lambda_j}\lbrack \kappa(2n-i-j)\rbrack_{\lambda_i+\lambda_j}}\\ &\quad\times\prod_{i=1}^n\frac{\lbrack a-1+\kappa(n-i)\rbrack_{\lambda_i}}{\lbrack a-1+2\kappa(n-i)\rbrack_{2\lambda_i}};
\end{split}
\end{equation*}
see Proposition 6.1 in Halln\"as. It is thus clear that the zeros, as well as the poles, of $Y_\lambda(0^n;a,\kappa)$ are located along hyperplanes given by $\kappa$ a negative rational number or $a -1 + \kappa i$ equal to a non-positive integer for some $i = 0,\ldots,2(n-1)$. This concludes the proof of Proposition \ref{rationalityProp}.

\subsection{Proof of Lemma \ref{convergenceLemma}}
In order to prove the statement we will make use of the so-called integral form
\begin{equation*}
	J_\lambda(x;\kappa) = \kappa^{-|\lambda|}h^\lambda P_\lambda(x;\kappa),\quad h^\lambda = \prod_{i=1}^{\ell(\lambda)}\prod_{j=1}^{\lambda_i}(\lambda_i-j+\kappa(\lambda^\prime_j-i+1)),
\end{equation*}
of the Jack polynomials. We recall that the coefficients in their expansion in the monomials $m_\lambda$ are polynomials in $1/\kappa$ with positive integer coefficients; see e.g.~Knop and Sahi \cite{KS97}. Since the Jack polynomials $J_\lambda$ are homogeneous of degree $|\lambda|$, it follows that
\begin{equation*}
	\left|J_\lambda(x;\kappa)\right|\leq J_\lambda(1^n;|\kappa|)||x||^{|\lambda|}
\end{equation*}
with $||x|| = \max(|x_i|)_{i=1}^n$. We recall from Stanley \cite{Sta89} (see Theorem 5.4) that
\begin{equation*}
	J_\lambda(1^n;\kappa) = \kappa^{-|\lambda|}\prod_{i=1}^{\ell(\lambda)}\prod_{j=1}^{\lambda_i}(j-1+\kappa(n-i+1)).
\end{equation*}
We observe that since $\kappa$ is not a negative rational number, there exist $\epsilon>0$ such that $|\kappa i + m|\geq\epsilon$ for all $i=1,\ldots,n-1$ and positive integers $m$. We fix such a number $\epsilon$. It is then clear that
\begin{equation*}
\begin{split}
	\prod_{i<j}\left|\frac{\lbrack\kappa(j-i-1)+1\rbrack_{\lambda_i-\lambda_j}}{\lbrack\kappa(j-i)+1\rbrack_{\lambda_i-\lambda_j}}\right| &= \prod_{i<j}\prod_{m_{ij}=1}^{\lambda_i-\lambda_j}\left|1 - \frac{\kappa}{\kappa(j-i)+m_{ij}}\right|\\ &\leq \prod_{i<j}\left(1+\frac{|\kappa|}{\epsilon}\right)^{\lambda_i-\lambda_j}\\ &\leq R_1^{|\lambda|}
\end{split}
\end{equation*}
for some $R_1$ independent of $\lambda$. Proceeding in a similar fashion, and using the fact that $J(1^n;|\kappa|)\leq |\kappa|^{-|\lambda|}\prod_{i=1}^{\ell(\lambda)}\lbrack |\kappa|n\rbrack_{\lambda_i}$, it is readily verified that also
\begin{equation*}
	\frac{|\kappa|^{|\lambda|}J_\lambda(1^n;|\kappa|)}{|h^\lambda|}\leq R_2^{|\lambda|}
\end{equation*}
for some $R_2$ independent of $\lambda$. We thus conclude that for some constants $C$ and $R$, independent of $\lambda$, we have
\begin{equation*}
\begin{split}
	|\mathscr{W}(x)| &\leq C\sum_\lambda \left|\prod_{i=1}^n\frac{\Gamma(a+\kappa(2n-i-1))}{\Gamma(a-1+\kappa(2n-i-1)+\lambda_i)}\right|(R||x^{-1}||)^{|\lambda|}\\ &\leq C \prod_{i=1}^n\sum_{\lambda_i=0}^\infty \frac{|\Gamma(a+\kappa(2n-i-1))|}{|\Gamma(a-1+\kappa(2n-i-1)+\lambda_i)|}(R||x^{-1}||)^{\lambda_i},
\end{split}
\end{equation*}
which clearly converges uniformly and absolutely on $\mathbb{C}\setminus N$. The statement thus follows from the so-called Weierstrass M-test; see e.g.~Theorem 11.5 in Apostol \cite{Apo67}.

\subsection{Proof of Theorem \ref{mainThm}}\label{proofThm21}
Let $\nu = (\nu_1,\ldots,\nu_n)$ be a partition. Using the orthogonality of the Jack polynomials, as well as the fact that $x_1\cdots x_n P_\nu(x) = P_{\nu+(1^n)}$ and $P_\nu$ is homogeneous of degree $|\nu|$, it is a straightforward exercise to verify that
\begin{multline}\label{F1Moments}
	\int_{T^n}P_\nu(x)W(x)dx\\ = (-2)^{|\nu|+n}\prod_{i<j}\frac{\Gamma(\kappa(j-i+1)+\nu_i-\nu_j)}{\Gamma(\kappa(j-i)+\nu_i-\nu_j)}\prod_{i=1}^n \frac{\Gamma(\kappa i)\Gamma(a+\kappa(2n-i-1))}{\Gamma(\kappa)\Gamma(a+\kappa(2n-i-1)+\nu_i)}.
\end{multline}
Since $\Gamma(\alpha+m)/\Gamma(m) = \lbrack\alpha\rbrack_m$, it follows from \eqref{F1Moments} and Proposition \ref{rationalityProp} that
\begin{equation}\label{F1BesselProd}
	\int_{T^n}Y_\mu(x)Y_\mu(x)W(x)dx = \Bigg(\prod_{i<j}\frac{\Gamma(\kappa)\Gamma(\kappa(j-i+1))}{\Gamma(\kappa i)\Gamma(\kappa(j-i))}\Bigg) f(a,\kappa)
\end{equation}
for some rational function $f$ with poles only at the excluded values of $a$ and $\kappa$. In order to verify that the right-hand side of \eqref{F1BesselProd} is zero, i.e., to establish the orthogonality relation \eqref{orthogonalityRel}, it is thus sufficient to verify that this function $f$ vanish on a subset of $\mathbb{C}^2$ which is dense in the Zariski topology.\footnote{We recall that a subset $Y\subset\mathbb{C}^2$ is dense in the Zariski topology if and only if the zero polynomial is the only polynomial which vanish on $Y$.} To this end we fix a non-negative integer $m$ such that $m\geq\max( |\lambda|,|\mu|)$, and consider the set of parameter values $(a,\kappa)\in\mathbb{C}^2$ such that $\kappa\geq 0$ and $a<-2(m+\kappa(n-1))+1$. It is clear that this subset of $\mathbb{C}^2$ is dense in the Zariski topology. We recall from Theorem 7.1 in Halln\"as \cite{Hal08} the following orthogonality relation for the multivariable Bessel polynomials of degree at most $m$: if $\lambda$ and $\mu$ are two partitions such that $|\lambda|,|\mu|\leq m$ and $\lambda\neq\mu$ then
\begin{equation}\label{L2OrthogonalityRel}
	\int_{\mathbb{R}_+^n}Y_\lambda(x)Y_\mu(x)W_{L^2}(x)dx = 0,\quad W_{L^2}(x) = \frac{1}{n!}\prod_{i=1}^n x_i^{a-2}e^{-2/x_i}\prod_{i<j}|x_i - x_j|^{2\kappa}.
\end{equation}
We stress that if either of the partitions $\lambda$ and $\mu$ has a weight greater than $m$ then this integral might not exist; see Lemma 7.1 in Halln\"as. We proceed to compute the integral $\int_{\mathbb{R}_+^n}P_\nu(x)W_{L^2}(x)dx$ for the partitions $\nu = (\nu_1,\ldots,\nu_n)$ such that $|\nu|\leq 2m$, and to show that, up to a non-zero constant, it is also given by the right-hand side of \eqref{F1Moments}. As we will show below, the value of this integral is easily inferred from a particular limiting case of the following generalisation of Selberg's integral formula due to Kadell \cite{Kad97}:

\begin{theorem}[Kadell]
Let $\nu = (\nu_1,\ldots,\nu_n)$ be a partition and let $\alpha,\beta,\kappa\in\mathbb{C}$ be such that
\begin{equation*}
	\Re(\alpha) > - \nu_n,\quad \Re(\beta) > 0,\quad \Re(\kappa)\geq 0.
\end{equation*}
Then
\begin{multline}\label{KadellValue}
	\frac{1}{n!}\int_{\lbrack 0,1\rbrack^n} P_\nu(y;\kappa)\prod_{i=1}^n y_i^{\alpha-1} (1 - y_i)^{\beta-1}\prod_{i<j}|y_i-y_j|^{2\kappa}dy\\ = \prod_{i<j}\frac{\Gamma(\kappa(j-i+1)+\nu_i-\nu_j)}{\Gamma(\kappa(j-i)+\nu_i-\nu_j)}\prod_{i=1}^n \frac{\Gamma(\alpha+\kappa(n-i)+\nu_i)\Gamma(\beta+\kappa(i-1))}{\Gamma(\alpha+\beta+\kappa(2n-i-1)+\nu_i)}.
\end{multline}
\end{theorem}

We recall that $(1-x/\beta)^\beta\rightarrow e^{-x}$ as $\beta\rightarrow\infty$, and also that
\begin{equation*}
	\lim_{\beta\rightarrow\infty}\beta^{b-a}\frac{\Gamma(\beta+a)}{\Gamma(\beta+b)}  = 1
\end{equation*}
for any complex $a$ and $b$; see e.g.~Theorem 3.4-1 in Carlson \cite{Car77} for the latter statement. By first replacing each variable $y_i$ by $y_i/\beta$, then using the fact that $P_\nu$ is homogenous of degree $|\nu|$, and finally taking the limit $\beta\rightarrow\infty$, it is thus straightforward to infer from \eqref{KadellValue} that
\begin{multline*}
	\frac{1}{n!}\int_{\mathbb{R}_+^n}P_\nu(y)\prod_{i=1}^n y_i^{\alpha-1}e^{-y_i}\prod_{i<j}|y_i-y_j|^{2\kappa}dy\\ = \prod_{i<j}\frac{\Gamma(\kappa(j-i+1)+\nu_i-\nu_j)}{\Gamma(\kappa(j-i)+\nu_i-\nu_j)}\prod_{i=1}^n\Gamma(\alpha+\kappa(n-i)+\nu_i);
\end{multline*}
c.f.~Section 17.3 in Mehta \cite{Meh91}. On the other hand, by a change of variables to $(y_1,\ldots,y_n) = (2/x_1,\ldots,2/x_n)$, it is easy to verify that
\begin{multline*}
	\int_{\mathbb{R}_+^n}P_\nu(x)W_{L^2}(x)dx\\ = \frac{2^{(a-1)n+\kappa n(n-1)+|\nu|}}{n!}\int_{\mathbb{R}_+^n}P_\nu(y^{-1};\kappa)\prod_{i=1}y_i^{-a-2\kappa(n-1)}e^{-y_i}\prod_{i<j}|y_i-y_j|^{2\kappa}dy.
\end{multline*}
We fix a non-negative integer $N$ such that $\nu\subseteq(N^n)$, and let $\hat{\nu}$ be the partition obtained by taking the complement of $\nu$ in $(N^n)$. Explicitly,
\begin{equation*}
	\hat{\nu}_i = N - \nu_{n-i+1},\quad i = 1,\ldots,n.
\end{equation*}
Then
\begin{equation*}
	(y_1\cdots y_n)^N P_\nu(1/y_1,\ldots,1/y_n) = P_{\hat{\nu}}(y_1,\ldots,y_n).
\end{equation*}
Indeed, the two sides have the same leading term, and it is readily verified that they satisfy the same orthogonality relation; see Macdonald \cite{Mac}. By combining the observations above it is now a straightforward exercise to establish the following:

\begin{proposition}\label{F2Prop}
Let $\nu = (\nu_1,\ldots,\nu_n)$ be a partition and $a,\kappa\in\mathbb{C}$ be such that
\begin{equation*}
	\Re(\kappa)\geq 0,\quad \Re(a) < -\nu_1 - 2\Re(\kappa)(n-1) + 1.
\end{equation*}
Then
\begin{multline*}
	\int_{\mathbb{R}_+^n}P_\nu(x)W_{L^2}(x)dx = 2^{(a-1)n+\kappa n(n-1)+|\nu|}\prod_{i<j}\frac{\Gamma(\kappa(j-i+1)+\nu_i-\nu_j)}{\Gamma(\kappa(j-i)+\nu_i-\nu_j)}\\ \times\prod_{i=1}^n \Gamma(-a+1-\kappa(2n-i-1)-\nu_i).
\end{multline*}
\end{proposition}

We observe the formula $\Gamma(-\alpha-m+1) = (-1)^m \Gamma(-\alpha+1)\Gamma(\alpha)/\Gamma(\alpha+m)$. By comparing Proposition \ref{F2Prop} with \eqref{F1Moments} we thus find that
\begin{multline}\label{integralEquality}
	\int_{T^n}P_\nu(x)W(x)dx\\ = \frac{(-1)^n}{2^{(a-2)n+\kappa n(n-1)}}\prod_{i=1}^n \frac{\Gamma(\kappa)}{\Gamma(\kappa i)\Gamma(-a+1-\kappa(2n-i-1))}\int_{\mathbb{R}_+^n}P_\nu(x)W_{L^2}(x)dx
\end{multline}
if the conditions on $a$ and $\kappa$ stated in Proposition \ref{F2Prop} are satisfied. We observe that the product of gamma functions preceding the integral in the right-hand side of \eqref{integralEquality} is non-zero. Since the Jack polynomials $P_\nu$ with $|\nu|\leq 2m$ span the space of symmetric polynomials of degree at most $2m$, the orthogonality relation \eqref{L2OrthogonalityRel} thus implies that the function $f$ in the right-hand side of \eqref{F1BesselProd} vanish for all $(a,\kappa)\in\mathbb{C}^2$ such that $\kappa\geq 0$ and $a<-2(m+\kappa(n-1))+1$, and hence is identically zero. This concludes the proof of the orthogonality relation \eqref{orthogonalityRel}. Moreover, by using \eqref{integralEquality}, and proceeding as above, it is now straightforward to infer the normalisation factor formula \eqref{normFactorFormula} from Lemmas 7.3 and 7.4 in Halln\"as \cite{Hal08}.

\section*{Acknowledgements}
I would like to thank J.~F.~van Diejen for an inspiring discussion on orthogonality, and A.~N.~Sergeev and A.~P.~Veselov for bringing Okounkov's binomial formula for the $BC_n$ Jacobi polynomials to my attention. I am also grateful to T.~H.Koornwinder for pointing out the references \cite{Dur93,EEKL93,KKH92}.

\bibliographystyle{amsalpha}

\end{document}